\input amstex
\def\Z{\text{\bf Z}}
\def\R{\text{\bf R}}
\def\N{\text{\bf N}}
\def\brcs#1{\left\{ #1 \right\}}
\def\Arrow #1;#2.#3{#1{}: #2 \to #3}
\def\st{such that }
\def\Set #1#2{\left\{#1 \vphantom{#2}\,\left|\vphantom{#1}\,#2\right. \right\}}







\font\rm=cmr10 \rm

\font\bf=cmb10
\font\Rm=cmr9 at 11pt
\rm
\font\it=cmsl9 at 10pt
 at 7pt

\font\Rrm=cmr17 at 16pt
   \font\Rm=cmr12 at 11.5pt

\long\def\Pf{\par\noindent {\it Proof.} }
\def\({\left(}
\def\){\right)}
\def\st{such that }
\def\qed{\hfill$\bullet$\vskip 4pt}

\def\brcs#1{\left\{ #1\right\}}

\def\wrt{with respect to }
\def\:{\,:}

\def\R{\text{\bf R}}
\def\N{\text{\bf N}}
\def\Z{\text{\bf Z}}

\def\Arrow #1;#2.{#1\:#2 \to }

\def\Set#1#2{\brcs{#1 \left|\vphantom{#1 #2} \right.#2}}

\def\Oh#1{{\pmb O}\(#1\)}

\def\oh#1{{\pmb o}\(#1\)}

\def\Rrr#1,#2{{\Cal J}_{#1,#2}}
\def\slfrac#1#2{{\raise -.07 ex\hbox{$^{#1}$}}\!/\raise .35 ex \hbox{${}_{#2}$}}
\def\ssf #1/#2{\slfrac {#1}{#2}}

\def\pd #1,#2.{\frac {\partial #1}{\partial #2}}

   \long\def\Lem
#1.#2\par{\vskip4pt{\baselineskip=13pt\font\it=cmsl12 at
11.5pt\Rm
   \noindent {\rm \uppercase{#1}} #2\vskip3pt

   }} 

\long\def\Proclaim #1.#2 \endproclaim{\vskip4pt{\baselineskip=13pt\font\it=cmsl12 at
11.5pt\Rm
   \noindent {\rm \uppercase{#1}} #2\vskip3pt

   }} 

\long\def\remark #1\endremark{\vskip 2pt \noindent {\it Remark\/} #1\par}

\long\def\Sectionhead #1.#2:\par #3{\vskip 4pt \noindent {\bf #1 #2}vskip 2pt\noindent\nospace #3}

\long\def\Title #1\par {\noindent{\Rrm #1}\vskip 9pt}

 \long\def\SubT #1.{\noindent {\it #1\/} } 
 
 \long\def\SecT
#1\par{\vskip 3pt \noindent {\bf #1}\vglue1pt
   \noindent}

\long\def\subtitle #1.{\vskip 2pt \noindent {\it #1}}

\long\def\Rmk#1\par{\vskip 1pt \noindent {\it
Remark.} #1\vskip2pt}

\long\def\Abstract #1\par{{\leftskip= 3 true cm \rightskip = 3 true cm \font\it=cmsl10 \font\rm=cmr10 \baselineskip = 10pt
\parindent=.35 true cm\rm\noindent 
{\it Abstract} #1\vskip 8pt

}}

\long\def\Author #1 \par{\noindent{\it #1}}
 

\def\flo #1{\lfloor #1 \rfloor}

\long\def\Ex #1.#2\par{\vskip 2pt \noindent {\it #1} #2}

\def\Gd{\subset {} \hskip -1.05em \raise .8pt \hbox{ \text{{\sm G}}}\,\ }

\NoBlackBoxes

\comment 
\long\def\text #1{{{\rm #1}}}
\def\frac #1#2{{#1 \over #2}}
\endcomment 

\def\oneone{1.1}
\def\onetwo{1.2}
\def\onethr{1.3}
\def\onefou{1.4}

\def\flo #1{\lfloor #1 \rfloor}

\def\tripnorm #1xxx{\left\|\hglue-.2ex\left|#1\right|\hglue-.2ex\right\|}

\def\flo #1{\lfloor #1 \rfloor}

\def\Oh#1{O\(#1\)}

\Title 
Random sets and intersections

\Abstract The following class of problems arose out of vain attempts to show that the Pascal's triangle adic transformation has trivial spectrum. Partition a set of size $N$ into sets of size $S \equiv S(N)$ (ignoring  leftovers). What is the likelihood that a set of size $K \equiv K(N)$ will intersect each set in the partition in at least $R \equiv R(N)$ members (as $N$ increases)? Via elementary techniques and under reasonable hypotheses, we obtain an easy-to-use formula. Although different from the corresponding minimum problem for balls and bins (with $m = K$ balls and $n = N/S$ bins), under modest constraints, the asymptotic probabiliities are the same.

\plainfootnote{}{\rm MSC 2010: 60C05 (primary); 05A16.}

\Author David Handelman %
\plainfootnote{$^{1}$}{\rm Supported in part by an NSERC Discovery Grant.}

\plainfootnote{${}$}{\hglue -2em  {\it Keywords and phrases\/}: random sets, inclusion/exclusion, Bonferroni's inequalities, log concave, strongly unimodal, balls and bins}

\vskip-6pt

\noindent Let $N$ be a (large) integer, and $S <N$, $R$, and $K$ positive integers. Chop the set $\brcs{1,2,3,\dots,N}$ into  $\flo{N/S}$ disjoint intervals of
length $S$ (that is, $\brcs{1,\dots , S}$, $\brcs{S+1, \dots, 2S}$, etc),
discarding the remainder. Regarding  $S \equiv S(N)$, $R\equiv R(N)$, and $K\equiv K(N)$  as (given) functions of $N$, we want to determine the likelihood that a subset of $\brcs{1,2,3,\dots,N}$ with $K$ elements will hit every one of the $\flo{N/S}$ subintervals in at least $R$ points (as $N\to \infty$). In order for the limiting probability to exist, there must be (mild) constraints on the variables, $S,R,K, N$. 

Here $G_j (t) = t^j e^{-t}$. 

\Lem Theorem \oneone. The limiting probability that a subset of $\brcs{1,2,\dots, N}$ of cardinality $K(N)$ will intersect each of the $\flo{N/S(N)}$ disjoint sets of cardinality $S(N)$ in at least $R(N)$ members is $ e^{-c}$, if 
$$c = \lim_{N\to \infty }\frac{N\cdot G_{R-1}\(\frac{SK}{N}\)}{S\cdot (R-1)!} 
$$
exists in $\R^+ \cup \brcs{\infty}$, and  provided that the following conditions hold.
\item{(a)} $R(N)^2 = \oh{ S(N)}$;
\item{(b)} $NR(N) = \oh{S(N)K(N)}$; 
\item{(c)} $R(N) S(N) \vee R(N) K(N) = \oh{N}$.

The method of proof is completely elementary, using inclusion/exclusion, Bonferroni's inequalities, and log concavity of various sequences.

\noindent {\it Examples.} It is convenient to pose the question in the form, given $R(N)$ and $S(N)$, how large must $K(N)$ be to guarantee occupation either  with almost certainty, or with probability $ 1/e$. 

The function  $G_{j}$, given by  $t \mapsto t^j e^{-t}$, if restricted to  $t \geq 
j$, is decreasing, and its image on this
interval is $(0,(j/e)^j)$.  We may define its inverse function, $\Arrow
T_j; (0,(j/e)^j). (j, \infty)$. This is a  relative of the Lambert
W-function. It is  decreasing, and on $(0,\epsilon)$, it expands
as $T_j(t) \sim -\ln t + j\ln |\ln t| + \Oh{\ln|\ln t|/|\ln t|}$.

Within the ranges permitted by  (a--c),  the probability of appropriately
hitting the all the subsets is monotone increasing in $K(N)$ and decreasing in $R(N)$ and $N/S(N)$. Beyond this range, monotonicity is obvious for $R(N)$  and $N/S(N)$, but not clear for $K$, since we are calculating the relative probabilities \wrt the number of $K$-element sets.

\Lem Examples. In these examples, we assume constraints (a--c) apply.

\noindent (a) If $R = 1$,   consider $(N/S) e^{-SK/N}$. To obtain $c = t$ (with $0 < t < \infty$), we would
have $K \sim (N/S)(\ln (N/S) - \ln t) $; to obtain $ c =
\infty$, we would have to have $KS/N -\ln (N/S)\to -\infty$ (in
addition to the constraint that $SK/N \to \infty$). In particular, if $S(N) = \sqrt N$, then $c = 1$ arises when $K(N)
\sim\frac 12 \sqrt N \ln N$, and $c = \infty$ (practical certainty of
unsuccessful occupation) occurs when 
$K/\sqrt N - \frac 12 \ln N \to -\infty$.

\noindent (b) When $R = r+1$ (constant in $N$), the situation is slightly more
complicated. For all sufficiently large $N$ (needed to ensure that $SK/NR$
is big), $K (N) \sim (N/S) T_r (t \cdot r! \cdot S/N)$. As $S/N \to 0$, we
can use the asymptotic expansion of $T_r$ (for very small $s> 0$), $T_r
(s) \sim -\ln s + r\ln|\ln s|$ (the next term is $\ln |\ln s|/ |\ln s| $, so is irrelevant); hence
$$\eqalign{
K(N) \sim \frac{N}S \(\ln \(\frac N{r!t S}\) + r \ln \ln \(\frac N{r!t
S}\)\)\cr
\sim \frac NS \(\ln \frac NS - \ln (r! t) + r \ln \ln \frac NS\).\cr
}$$
(The next term is around $rt!\ln\ln (N/S)/\ln (N/S)$.) This is amazingly sensitive to
changes in $t$, but at least $K$ is infinitely bigger than the $K$
required for $R = 1$. See also the {\it Perturbations and asymmetry\/} comment below.

\noindent (c) With unbounded $R(N)$, computations are naturally more complicated, but we
can see what happens in the special case that $S(N) = \sqrt N$. Set $r \equiv r(N)  =  R(N) +1$, and set $K = g(N)r\sqrt N $; in order for Theorem \oneone\ to apply, we require (at least) $g(N) \to \infty$. We will show that $ c = \infty$ (that is, the probability is $0$) if $r(N)g(N) \leq \frac 12 \ln N$, and $c > 0$ if  $r(N) g(N) \geq (\frac12 + \eta) \ln N$  for some $\eta > 0$. If $r(N) \geq k \ln N$ for some real $k> 0 $, we are out of luck (as $g(N) \to \infty$ is required for our result, in order that the hypothesis (b), $RN = \oh{SK} $, be satisfied). 

To compute $c$, using Sterling's formula, $c = \lim_N f(N)$, where
$$\eqalign{
f(N) &= \frac{\sqrt N (rg)^r \exp(-rg)}{\(\frac{r}e\)^{r}\sqrt {2\pi r}}, \qquad\text{so} \cr
\ln f(N) &=\frac 12 \ln N + r(\ln g - g) + r- \frac
12 \ln (2\pi r).\cr
}$$
As $g(N) \to \infty$, it follows that  $rg \sim r(g - \ln g)$. 

If $rg\leq \frac 12 \ln N$, then since $r(N)\to \infty$, we see that $\ln f(N) \to \infty$, so $c = \infty$. On the other hand, if $r(N)g (N) \geq (\frac 12 + \eta) \ln N$, then $\ln f(N) \leq - \eta \ln N/2  $; 
so $\ln f(N ) \to -\infty$, and thus $c = 0$. 
\qed

\noindent{\it Perturbations and asymmetry.} Asymmetric behaviour appears when we 
consider perturbations. Here, we permit $R(N)$ to be unbounded (but
still subject to all the constraints appearing in Theorem \oneone, e.g., $RN = \oh{SK}$). For all
sufficiently large $N$, define $K^0 (N)$ to be the solution to $c = 1$,
that is (with $r(N) = R(N)-1$),
$$
\frac{(N/S) (SK^0/N)^r \exp (-SK^0/N)}{r!} = 1.
$$
Now we perturb $K^0$, defining $K(N) = K^0(N) + a N/S(N)$ where $a$ is a real number
independent of $N$; it can be of either sign. We calculate the new value
for $c$,
$$ \eqalign{
c_1 (N)& = \frac{(N/S) (SK^0/N + a)^r \exp (-SK^0/N +a)}{r!}\cr
& = 1\cdot\(1 + \frac {aN}{SK^0} \)^r e^{-a}\cr
& \sim \exp \(raN/SK^0 \) \cdot e^{-a} \cr
& \to e^{-a}. \cr
}$$

Next, we permit $a$ to vary. If $a \to \infty$, then $ \exp\({-c_1}\) = 
e^{-e^{-a}} \sim 1- e^{-a}$---in other words, convergence (in $a$) entails
exponential convergence (in $a$) of the probability to $1$. On the other hand, if
$a \to -\infty$, then $e^{-c_1} = \exp(-e^{|a|})$. This time, we have
{\it ultra-\/}exponential decay (in $a$) to $0$. So there is asymmetry in the rates
of convergence to one and to zero.\qed

\noindent {\it Significance of conditions a--c.\/} Condition (b) seems very reasonable in the sense that it appears to be close to necessary for the probability to be nonzero (assuming that $K = \oh{N}$, which is itself reasonable); however, I have not been able to rule out the possibility that $\alpha(N)$ be bounded (but large),  $R(N)\to \infty$, $K = \oh N$, and the likelihood of successful occupation being nonzero. The constraints in (a,c) turned up in the course of the proof. It would be interesting to clean up these extra conditions. \qed

\comment
Suppose $R(N) = 1$ is constant in $N$. Then $c = (N/S)e^{-SK/N}$; if $c$ is finite, then $K= (\ln c + \ln (N/S))(N/S)$, whereas if $ c = \infty$, $K(N) \gg (N/S) \ln (N/S)$. A natural case occurs when $S = \lambda \sqrt N$ for some positive constant $\lambda$. Then $K(N) \sim \lambda^{-1}\sqrt N ((\log N)/2 - \ln \lambda)$ will guarantee $c =1$, that is the likelihood that a random set of this size  will hit all the intervals tends to $1-1/e$. And $c = \infty$ corresponds to $K(N) - (\lambda^{-1}\sqrt N ((\log N)/2 - \ln \lambda)) \to \infty$. 

When $R(N) = r +1 > 1$ is constant in $N$, and say $c = 1$ (so the probability will tend to $1-1/e$), we have $K(N) \sim  (N/S)T_{r}((S/N)r!)$. Since $S = \oh{N}$, the asymptotic expansion of $T_r$  is valid, and so the right side is about $(N/S) (\ln (N/S) - \ln r! + r\ln\ln (N/S))$. This is bigger than the preceding, but not by very much. 

When $R(N) \to \infty$, the situation is more complicated, and not just because there are constraints on the growth. For example, suppose $R(N) \sim \ln N$. For the case $c = 1$, we have $K = (N/S)T_{R-1}((R-1)! S/N)$. For the asymptotics of $T$ to be valid, we require $(R-1)! S/N$ to be small, which is not going to happen (since $\Gamma (\ln N)$ is a lot bigger than $N$). Fortunately, $T_j$ is defined on $(0,(j/e)^j)$, and when $ j \sim \ln N$, $(j/e)^j \sim \Gamma(\ln N +1)$, and so $(R-1)! S/N = \oh{(j/e)^j}$; that is, it lies well within the domain of $T_{R-1}$. I don't at the moment see how to estimate it. 

More examples?

\noindent {\it Asymmetry.} Suppose that $R(N)$ and $S(N)$ are given (for all $N$), and $K^1(N)$ is defined so that the corresponding value of $c$ is $ 1$. Consider what happens when we replace $K^1$ by $K^1 + a$ (where $a \equiv a(N)$ can either be positive or negative); let $c_a$ be the corresponding value (if it exists at all) with this new choice for $K$. 

If $R > 1$, the derivative of $G_{R-1}$ is $t\mapsto t^{R-2}((R-1) - t) e^{-t}$, so if we view $c$ as a function of $K$, $\frac{dc}{dK} (SK/N) = (1- SK/N(R-1))G_{R-2} (SK/N)/(R-2)! \sim -(SK/NR)G_{R-2}(SK/N)/(R-2)!$.  aargh

\noindent {\it Hypotheses.} Hypotheses (a) and (c) look reasonable, in view of the original problem. Hypotheses (b) and (d) are there so that the expression simplifies (see the proof). Without them, the expression for $c$ would become tremendously awkward.

\item¥
discussion of hypotheses

\item¥ asymmetry around critical point

\item¥ connections to bins and balls

\item¥ proof of corresponding cases for bins and balls

\endcomment 

We begin the proof of Theorem \oneone.

For each $i$, we  define  $\beta_i(N)$ for the purposes of using inclusion/exclusion, as follows. Pick a collection of $i$ of the $N/S$ intervals, and determine the number of subsets of size $K \equiv K(N)$ that hits each interval in the collection in fewer than $R \equiv R(N)$ elements. Then multiply by ${{N/S} \choose i}$ (the number of choices of $i$ of these intervals. Finally divide by the number of $K$-element subsets of ${1,2,\dots, N}$, that is, divide by ${N \choose K}$.  Then $1-\sum_{i=1}^{\infty}(-1)^{i+1} \beta_i(N)$ (if it converges) is the
probability, obtained from inclusion/exclusion, that a set of size $K$ will
hit all of the disjoint intervals of size $S(N)$, except possibly the last, short, 
one (if $S$ does not divide $N$), in at least $R$ points.

As a consequence of Bonferroni's inequalities, there is a unexpected but  pleasant convergence property.

\Lem Lemma \onetwo. Suppose that $\beta_i(N) \to b_i $ for all $i$, and that
$\sum_{i=1}^{\infty} (-1)^{i+1} b_i$ converges to $C$. Then
$\sum_{i=1}^{\infty} (-1)^{i+1} \beta_i(N) $ converges (as $N \to \infty$) to $C$.

\Pf Set $f(N)= \sum_{i=1}^{\infty} (-1)^{i+1}\beta_i(N) $ (for each $N$, this is a
finite sum). For all $k,l$, by Bonferroni's inequalities, we have
$$
\sum_{i=1}^{2k} (-1)^{i+1} \beta_i (N) \leq f(N) \leq \sum_{i=`1}^{2l+1}
(-1)^{i+1} \beta_i (N).
$$
For each $k,l$, the sums on the left and right converge as $N \to \infty$
to $\sum_{i=`1}^{2k} (-1)^{i+1} b_i$ and $\sum_{i=1}^{2l+1} (-1)^{i+1}
b_i$, respectively. On taking the $\liminf $ on the left and the $\limsup$ on the right, we have that  for all $k$ and $l$,
$$
\sum_{i=`1}^{2k} (-1)^{i+1} b_i \leq \liminf_N f(N) \leq \limsup_N f(N)
\leq \sum_{i=`1}^{2l+1} (-1)^{i+1} b_i.
$$
By hypothesis, the left and right sides converge as $k,l \to \infty$ to
the same thing, so $\liminf f(N) = \limsup f(N)$.
\qed

This simplifies the arguments, as it is relatively easy to calculate $b_m:=
\lim_{N \to \infty}\beta_m(N)$.

\def\ch#1,#2 {{{#1}\choose {#2}}}

Fix $N$, $K = K(N)$, $S = S(N)$ and $R = R(N)$, and $m$. Varying
over intervals of the form $[kS, (k+1)S) \cap \N$, we have an exact, but
horrendous, formula for the $m$th term arising from inclusion/exclusion.
For each $s = 0,1, 2, \dots, m\cdot (R-1)-1$, set
$$
D_R(s) \equiv D(s) = \Set{(i_j) \in \Z_+^m}{\sum i_j = s; \text{and  $i_j \leq
R-1$ for all $j$, }}.
$$
Then
$$
\beta_m(N) = \frac{\ch N/S,m \sum_{s=0}^{m\cdot(R-1)}\ch N-mS,K-s \sum_{(i_j) \in D(s)}
\prod_{j=1}^m \ch S,i_j }
{\ch N,K }.
$$
The horrible sum will be considerably simplified---under modest
conditions, we can delete all the terms $s = 0, 1, \dots, m(R-1) -1$, leaving
just one term arising from $s = m\cdot (R-1)$, that is, $\(\ch S,R-1 \)^m \ch
N-mS,K-m(R-1) $. This permits us to compute $b_m = \lim \beta_m(N)$
easily.

We wish to obtain  reasonable approximations for $\beta_m(N)$. We assume $R(N)^2 = \oh{K(N)}$ (this, and similar expressions, will usually be abbreviated by dropping the $N$---thus, $R^2 = \oh{K\wedge S}$), $NR = \oh{SK}$, $RS \vee RK = \oh{N}$, and $S \vee K = \oh{N}$.  

If $f$ is a polynomial, we use inner product notation, $(f,x^s)$, to denote the coefficient of $x^s$ appearing in $f$. Let $p$ the initial segment of degree $R-1$ of $(1+x)^S$. Then $(p^m, x^s)
= \sum_{(i_j) \in D(s)} \prod_{j=1}^m {S \choose {i_j}}$. Since $(1+x)^S$
is strongly unimodal (log concave), so is $p$ and thus so is $p^m$. In
particular, the sequence $s \mapsto \sum_{(i_j) \in D(s)} \prod_{j=1}^m {S
\choose {i_j}}$ is strongly unimodal. So is the sequence $s \mapsto {{N-mS}
\choose {K-s}}$. Since the Hadamard product of strongly unimodal sequences
is itself strongly unimodal, so is
$$
s \mapsto {{N-mS} \choose {K-s}}\cdot \sum_{(i_j) \in D(s)} \prod_{j=1}^m {S
\choose {i_j}}: = g_N(s). \tag 1
$$
This is over the range $s = 0, 1, \dots, (R-1)m$. In particular, if
$g((R-1)m)/g((R-1)m -1) = t >1$, then $g((R-1)m)< \sum_{0}^{(R-1)m} g(s)
\leq g((R-1)m)/(1-1/t)$ (since the sequence $(g(s))_{s=0}^{(R-1)\cdot m}$ is strongly unimodal).

Define $\alpha(N) = SK/RN$. When $s = (R-1)\cdot m$, $D(s)$ consists of the one point, $(R-1,R-1,\dots, R-1)$, and when $s = (R-1)\cdot m-1$, $D(S)$ consists of  $m$ points, $\brcs{(R-1,\dots, R-1, R-2,R-1, \dots, R-1)}$. Thus it is easy to calculate
$$\eqalign{
\frac {g((R-1)m)}{g((R-1)m -1)}& = \frac{\({S \choose {R-1}}\)^m
{{N-mS}\choose{K-(R-1)m}}}
{\({S \choose {R-1}}\)^{m-1}m {S \choose {R-2}} {{N-mS}\choose{K-(R-1)m
+1}}}\cr
& = \frac{(S-R+2)\cdot \(K-(R-1)m+1\)}{(R-1)m\cdot \( N-mS -K +(R-1)m\)} \cr
& \geq \frac{SK}{NRm} \(1 - m\,\Oh{\frac{R}K} \) \cr
& = \frac{\alpha(N)}m \( 1 - m\,\Oh{\frac{R}K}\).\cr
}$$
From (b), $R/K= \oh{S/N}$; from (c), $S = \oh{N}$, and it follows that $R = \oh{K}$.
Since $m$ is fixed and $\alpha (N) \to \infty$ by (b), the expression   becomes arbitrarily
large. Thus given $\epsilon$, there exists $N_0$ \st for all $N \geq
N_0$, $\sum g(s)/g((R-1)m) < 1+ \epsilon$, that is, almost all the mass of the sum in the definition of $\beta_m (N)$ is concentrated on the last term. [If $R=1$ or $2$, the
corresponding argument is trivial.]

Thus
$\beta_m(N)$ is asymptotic (in $N$) with
$$
\beta'_m(N):= \frac{\ch N/S,m \(\ch S,R-1 \)^m \ch N-mS,K-m(R-1) }{\ch N,K }.
$$

Now we can estimate $\beta'_m$. The following is elementary.

\Lem Lemma \onethr. If $B < A/2$, then
$$\eqalign{
\frac{A!}{(A-B)!} & = A^B \exp\(-\frac{B^2 - B}{2A} \)\(1 -
\Oh{\frac{B^3}{A^2}} \)\cr
{A\choose {B} } &= \frac{A^B}{B!} \exp\(-\frac{B^2 - B}{2A} \)\(1 -
\Oh{\frac{B^3}{A^2}} \).\cr
}$$

\Pf Obviously $A!/(A-B)! = \prod_{i=0}^{B-1}(A-i) = A^B
\prod_{i=0}^{B-1}(1-i/A)$. Now take logarithms of the product, and sum.
\qed

The constant(s) that appear in the big Oh terms are  slightly bigger than one-half. Obviously, finer estimates can be made, but we won't require them in what follows.

 Ignoring some of the $1 \pm \Oh{\cdot}$
terms,
$$\eqalign{
\beta'_m(N) & = \frac{(N/S)^m S^{m(R-1)}
}{m! ((R-1)!)^m} \exp\(-m^2\frac{(R-1)^2}{2S}\) \frac{\prod_{i=0}^{K-1}
(K-i)}{\prod_{i=0}^{K-1} (N-i)} \cdot \prod_{i=1}^{K+m(S-R+1)} \frac{N-mS -i}{K-m(R-1)-i}\cr
& = \frac{(N/S)^m S^{m(R-1)}
}{m! ((R-1)!)^m} \(1-m^2\Oh{\frac{R^2}{S}}\) \frac{\prod_{i=0}^{m(R-1)-1}
(K-i)}{\prod_{i=0}^{mS -1} (N-i)} \cdot \prod_{i=K}^{K+m(S-R+1)} (N -i)\cr
& = \frac{(N/S)^m S^{m(R-1)}
}{m! ((R-1)!)^m} K^{m(R-1)} N^{-mS}N^{m(S-R+1)} \times
\cr
& \qquad \exp\(-m^2\frac{(R-1)^2}{2K} + m^2
\frac{S^2}{2N} - m(S-R+1)\frac KN - m^2\frac{(S-R+1)^2}{2N}\) (1+ \oh{1}) \cr
& = \frac{\((N/S)\cdot (KS/N)^{R-1}e^{-SK/N}/(R-1)!\)^m}{m!}  \times
\cr
& \qquad \exp\(-m^2\frac{(R-1)^2}{2K} + m^2\frac{2SR-R^2-2S+2R-1 }{2N}
- m(S-R+1)\frac KN \) (1+ \oh{1}) \cr
&= \frac{\((N/S)\cdot (KS/N)^{R-1}e^{-SK/N}/(R-1)!\)^m}{m!}  \times
\cr
& \qquad \exp
\(- m\frac {SK}N \)\(1+m^2\Oh{\frac{R^2}{K} \vee \frac{SR}N}\) \(1 + m\Oh{\frac{RK}N}\) (1+ \oh{1}) \cr
& =\frac{\((N/S)\cdot (KS/N)^{R-1}e^{-SK/N}/(R-1)!\)^m}{m!} \(1 \pm
\oh{1}\)\cr
}$$
To go from the penultimate line to the last line, we observe that $NR = \oh{SK}$  entails $R^2 = \oh{SRK/N} = \oh{K}$ (from (c)). Thus,
$$\eqalign{
\beta'_m (N) &= \frac{\((N/S)\cdot 
(R\alpha(N))^{R-1}e^{-R\alpha(N)}/(R-1)!\)^m}{m!} \(1 \pm \oh{1}\)\cr
& = \frac{\(\frac{N\cdot G_{R-1}\(SK/N\)}{S\cdot (R-1)!}\)^m}{m!} \(1 \pm
\oh{1}\),\cr
}$$
where $G_{j} (x) $ is function $t \mapsto t^j e^{-t}$, defined for $t >
j$. 

Hence, if the sequence $((N/S)\cdot G_{R-1}(SK/N)/(R-1)! ) 
$ converges to some number $c$, then $b_m = c^m/m!$, and by Lemma \onetwo,
$\sum_{i=1}^\infty (-1)^{i+1}\beta_i (N)$ converges to $ e^{-c}$. In
this case, the likelihood that a set with cardinality $K(N)$ hits each of
the $S(N)$ equi-length intervals in at least $R(N)$ points converges to
$e^{-c}$.

If $((N/S)\cdot G_{R-1}(SK/N)/(R-1)! ) \to \infty$, a
routine argument shows the limiting probability is zero, as follows.
First, assume $R(N) > 1$ for almost all $N$. We may select an increasing
sequence of positive real numbers, $t_n \to \infty$ with the property that
for all $N \geq n$, $t_n \sqrt{2\pi (R(N)-1)}< N/S(N)$ (possible, since $RS= \oh{N}$). For $N \geq n$, consider $c(N,n): = t_n(R(N)-1)!S(N)/N $. We have
$c(N,n) < (R(N)-1)!/ \sqrt{2\pi (R(N)-1)} \leq ((R(N)-1)/e)^{R(N)-1}$.
Hence $T_{R(N)-1} (c_{N,n})$ is defined.

For $N \geq n$, set $\alpha_{n}(N) = T_{R(N)-1} (c_{N,n})/R(N)$. Now fix
$n$, and consider the sequence 
$$\( \frac{N\cdot G_{R(N)-1} (\alpha_n
(N)\cdot R(N))}{S(N)\cdot (R(N)-1)!}\)_{N \geq n}.
$$
 By construction, this converges to $t_n$. On
setting $K_n (N) = R(N)N\alpha_n(N)/S(N)$, we obtain a new set of
parameters, but only the $K$ terms are changed, and the likelihood of
successful occupation is $e^{-t_n}$. Since $T_{R(N)-1} $ is decreasing on
its ray of definition, and $(t_n)$ is increasing, we see that $\alpha_n
(N) \geq \alpha_{n+1}(N)$ (when $N \geq n$), and so $K_n \geq K_{n+1} $.
In particular, $K \leq K_{n}$ for all $n$. Hence the likelihood of
successful occupation for the original sequence of parameters $(N,S(N),
K(N), R(N))$ is bounded above by  the corresponding likelihoods for
each of $(N,S(N), K_n (N), R(N))$. Since the latter likelihood is
$e^{-t_n}$ and this converges to zero, the likelihood of successful
occupation for the original set of parameters is zero.

If $R (N) = 1$ for all $N$, then the argument is even simpler, since we
only have to consider the inverse function of $x \mapsto e^{-x}$.

\noindent {\it Ball and bins.}
Unsurprisingly, the questions just considered are close to, but different
from, balls and bins problems. The closest relative, in its simplest form,
has $n$ bins and $m$ balls; the balls are tossed randomly into the
bins. What is the probability that the minimum number of balls in each
bin is at least $R$? We will replace $m$ by $K$. With $N/S = n$, this is
different from the preceding problem, in that in the former case, the
likelihood that a ball will land in one of the bins (one of the intervals,
$[kS, (k+1)S)$) depends on how many are already there. However, it turns
out that the asymptotic probabilities are identical.

Balls and bins has a large literature, mostly concerned with optimal
strategies,
or with the expected maximum number of balls. The basic reference is [RS], although Google Scholar lists well over one thousand papers containing the term {\it balls and bins\/} (of course, not all are relevant). I could not find a reference
to the expected minimum number of balls, presumably because it is not of interest.

We can easily apply the methods here to the minimal balls and bins
problem. Replace $m$ by $K$, and eventually, $n$ by $N/S$ (to conform with
the notation for sets). Suppose we have a set of $l$ bins (corresponding to $m$ in the former
context, as in $\beta_m$, but it would be too confusing to maintain this
convention).

There are $m = K$ balls in $n = S/N$ bins; what is the probability that all bins
have at least $r+1 = R$ balls?

Fix $l \leq n$. We obtain an expression for the probability that each the
first $l$ bins contain at most $r$ balls, and use inclusion/exclusion as in the previous situation.

Let $(a(1), \dots, a(l); b(l+1), \dots, b(n))$ be a distribution of the
$m$ balls (so that $\sum a(i) + \sum b(j) = K$). Then the likelihood of
obtaining this distribution of balls is exactly
$$
\frac{ {{K} \choose {a(1), \dots, a(l); b(l+1), \dots, b(n)}}}{n^K}. \tag *
$$
(The total number of distributions, the denominator, is of course given by
the sum of the coefficients of $(\sum_{i=1}^n x_i)^K$.) Temporarily fix
$a = ( a(i))$, the initial segment. We want a formula for the likelihood
that $a$ appears as the initial segment of the distribution, that is, the
sum over all admissible choices of sequences $(b(l+1), \dots, b(n))$ of
the expression appearing in (*) (with $a$ fixed).

Set $s = \sum a(i)$; then $\sum b(j) = K-s$; let $B $ be the set of
ordered $n-l$-tuples of nonnegative integers, $b = (b(l+1), \dots, b(n))$
whose sum is $K-s$. Then
$$\eqalign{
\sum_{b \in B} {{K} \choose {a(1), \dots, a(l); b(l+1), \dots, b(n)}}&=
\sum_{b \in B} \frac{K!}{\prod a(i)! \prod b(j)! }\cr
& = \frac{K!}{\prod a(i)!(K-s)!} \sum_{b \in B} \frac{(K-s)!}{\prod b(j)!}
\cr
& = \frac{K!}{\prod a(i)!(K-s)!} (n-l)^{K-s} \cr
}$$

Hence the probability of obtaining the initial sequence $a = (a(i))$  is (with $s =
\sum a(i)$),
$$
p(a):= \frac{K! (n-l)^{K-s}}{(K-s)! n^K} \cdot \frac 1{\prod a(j)!}.
$$
Set
$$D(s) = \Set{a \in \Z^l_+}{a(i) \leq r \text{ for all $ i = 1, 2, \dots,
l$, and $\sum a(i) = s$}}.
$$
Thus
$$
q(s):= \sum_{a \in D(s)} p(a) = \frac{K!}{\(\frac{n}{n-l}\)^K} \cdot
(n-l)^{-s} \cdot \frac 1{(K-s)! } \cdot \sum_{a \in D(s)} \frac 1{\prod
a(j)!}.
$$
The functions $s \mapsto (n-l)^{-s}$ and $s \mapsto 1/(K-s)!$ are log
concave. Since $\sum_{a \in D(s)} 1/\prod a(i)!= \( \(\sum_{t=0}^r
x^t/t!\)^l,x^s\)$, the function $s \mapsto \sum_{a \in D(s)} 1/\prod
a(i)!$ is also log concave. As the Hadamard product of log concave
functions defined on the integers is itself log concave, we obtain that $s
\mapsto q(s)$ is log concave.

Now define $\beta_l \equiv \beta_l (N)$ to be the sum over all sets of
$l$ bins, of the likelihood that particular set has at most $r$ balls in
each bin. This is obviously
$\beta_l = \sum_{s=0}^{rl} {n \choose l} q(s) = {n \choose l}
\sum_{s=0}^{rl} q(s)$.

If $q(rl)/q(rl-1) $ exceeds $t > 1$, then from strong unimodularity, we
see that $q(rl) \leq \sum q(s) \leq q(rl)/(1-1/t)$. Hence $\sum q(s) =
q(rl)(1+ \Oh{1/t}$. It is easy to calculate $q(rl)$ and $q(rl-1)$, since
$D(rl) = \brcs{(r,r,\dots, r)}$ and $D(ql-1)$ consists of the $l$ points
$\brcs{(r, r, \dots, r, r-1,r, \dots, r)}$. Thus
$$\eqalign{
\frac{q(rl)}{q(rl-1)} &= \frac{(n-l)^{rl-1}}{(n-l)^{rl}} \cdot
\frac{(K-rl+1)!}{(K-rl)!} \cdot \frac{l ((r-1)!)^{l-1}r!}{(r!)^l}\cr
& = \frac{K-rl+1}{n-l} \cdot \frac{l }{r}.\cr
}$$
We thus assume that (holding $l$ fixed) $K/nr \to \infty$ (this
corresponds to $KN/SR \to \infty$). (If $r = 0$, a separate, but easier,
argument is required.)

Thus $\sum_{s} q(s) = q(rl)(1+\Oh{nr/K})$, and so (up to multiplication by
$1+\Oh{nr/K}$ and with $l$ fixed),
$$\eqalign{
\beta_l &= \(\frac{n-l}{n}\)^K \cdot {\prod_{j=0}^{rl-1} (K-j) } \cdot {n
\choose l} \cdot \frac 1{(r!)^l} \cdot \frac 1{(n-l)^{rl}}\cr
& = \frac{ \exp (-Kl/n - Kl^2/2n^2) K^{rl} \exp(-(rl)^2/2K)) n^l}{l!
(n-l)^{rl} (r!)^l} \( 1 \pm l^3\Oh{\max\brcs{\frac K{n^3},
\frac{r^3}{K^2}, \frac rn}}\) \cr
&= \frac{\(\frac{n\(\frac Kn\)^r \exp-\frac Kn}{r!} \)^l}{l!}\(1 + l^3 \Oh{\max\brcs{\frac K{n^2},
\frac{r^2}{K}, \frac{r^2}n}}\).\cr
}$$
Thus, assuming $r^2 = \oh{K \vee n}$, $K = \oh{n^2}$, $nr = \oh{K}$, the formula is exactly the same as in the previous situation (with $n = N/S$). If $\lim n G_r(K/n)/r! = c \in \R^+ \cup \brcs{\infty}$, then the asymptotic probability is $e^{-c}$.
    
Properly stated, we have the following. Recall that $G_r(x) = x^r e^{-x}$.

\Lem Theorem \onefou.  Suppose that at discrete time $t$, $m \equiv m(t)$ balls are thrown at $n \equiv n(t)$ bins. Suppose that $r\equiv r(t)$ is a positive integer-valued function \st 
$$
\lim_{t \to\infty}  \frac{n G_r(m/n)}{r!} = c \in \R^+\cup \brcs{\infty}.
$$
Suppose in addition that $r = \oh{\sqrt n \vee \sqrt m}$, $m = \oh{n^2}$, and $nr = \oh{m}$.
Then the limit of the likelihoods that every bin contain at least $r(t) +1$ balls at time $t$, is $e^{-c}$.

\noindent {\it Motivation.} These result were motivated by an (unsuccessful) attempt to show that the Vershik adic map associated to Pascal's triangle is weakly mixing, that is, has trivial spectrum. The idea was that as a path progressed along the triangle, the chances that $e^{2\theta \pi i }$ belong to the spectrum would diminish. However, this depends on the behaviour of the sequence fractional parts, $\brcs{c(m)\theta}$, where $c(m) = {{2m} \choose m}$ is the central binomial coefficient, for arbitrary irrational $\theta$. This seems impossible to deal with.

\long\def\Rf[#1] #2, #3. #4\par%
{\parindent=10.5pt\vskip 4pt \itemitem{[#1]} #2, {\it #3,} #4\par\vskip2pt}

\SecT Reference

\Rf [RS] M Raab \& A Steger, Balls \& bins---a simple and tight analysis. Lecture Notes in Computer Science 1518 (1999) 159--170.

\vskip 8pt
\noindent Mathematics Department, University of Ottawa, Ottawa ON K1N 6N5 Canada;  dehsg\@uottawa.ca.

\vskip 4pt

\noindent Comments and criticisms are appreciated, especially if this turns out to have been done before!
\end

With $p = 1/n$ ($= S/N$), the likelihood of a bin
containing at most  $r = R-1$ balls is $\sum_{i=0}^{r} p^i (1-p)^{K-i},$
obtained as $\((1-p + px)^{K}, \sum_{i=0}^{R-1} x^i \)$, that is, the sum
over the tail $i \leq R-1$ of the coefficients in the binomial expansion.

Let $s \leq l\cdot (R-1)$, and define $D(s) $ as we did
previously. Define $q = lp/(1+lp)$. Then the likelihood each of this particular  set of bins has
fewer than $R-1$ balls and there are a total of $s$ balls in this set of bins, is
$$\eqalign{
A(l,s)&:= \(\(1-q + qx\)^K,x^s\)\times \frac{\sum_{a \in D(s)} {s \choose{a(1),a(2), \dots, a(l)}}}{l^s}\cr
& = q^s (1-q)^{K-s}{K\choose s}\frac{s!}{l^s} \cdot \sum_{a \in D(s)}\prod_{j=1}^l \frac{1}{ a(j)!}.\cr
}$$
The first line comes from the likelihood that exactly $s$ of the $K$ balls belong to the specified $l$ bins (the likelihood that a single ball so belongs is $lp$), multiplied by the number of acceptable configurations of the $s $ balls (the numerator the second term) divided by all possible configurations. Then $\beta_l (N) = {{N/S} \choose l}\cdot \sum_{s \leq l\cdot r}A(l,s)$. 

The functions $s\mapsto \(\(1-q + qx\)^K,x^s\)$, $s \mapsto s!/l^s$, and $s \mapsto \sum_{a \in D(s)}\prod_{j=1}^l \frac{1}{ a(j)!}$ are strongly unimodal. Strong unimodularity of the  first two is obvious, and that of the third follows from 
$$
\(\(\sum_{i=0}^r \frac{x^i}{i!}\)^l, x^s\) = \sum_{a \in D(s)}\prod_{j=1}^l \frac{1}{ a(j)!}.
$$
Thus the function $s\mapsto A(l,s)$ is the Hadamard product of strongly unimodal functions, so is itself strongly unimodal. Now
$$\eqalign{
\frac{A(l,rl)}{A(l,rl-1)}& = \frac{q^{rl} (1-q)^{K-rl}{K\choose rl}\frac{(rl)!}{l^{rl}\cdot (r!)^l} }{q^{rl-1} (1-q)^{K-rl+1}{K\choose rl-1}\frac{(rl-1)!l}{l^{rl-1}\cdot (r!)^{l-1}(r-1)!} }\cr
& = \frac q{1-q}\cdot  \frac{K-rl+1}{rl} \cdot \frac{rl}{l^2 r}\cr
& = \frac q{1-q}\cdot  \frac{K-rl+1}{rl^2}.\cr
}$$
With $l$ fixed, and $r =\oh{Kp}$ (corresponding to $NR = \oh{KS}$), we see that this blows up as $N \to \infty$. As in the preceding case, we deduce that $\sum_{s \leq lr} A(l,s) = A(l,rl) (1+ \Oh{r/Kp})$. 

Thus 
$$\eqalign{
\beta_l (N) &= {{N/S} \choose {l}}q^{rl} (1-q)^{K-rl}{K\choose rl}\frac{(rl)!}{l^{rl}\cdot (r!)^l} \(1+ \Oh{r/Kp}\)\cr
&= \frac{p^{-l} p^{rl} K^{rl} (rl)!}{l!(1+lp)^K (rl)!l^{rl}(r!)^l}\exp{}\ \(1 + \Oh{\max\brcs{\frac r{Kp},{l^3 p^2},\frac{r^3l^3}{K^2}}}\)
}$$

$$
A(l):= \sum_s \sum_{w = (a(j)) \in D(s)}  \(\(1-lp + p \sum_{i=1}^l
x_i\)^{K}, \sum x^w \),
$$
where we are using monomial notation, $x^w = x_1^{w(1)}\dots x_l ^{w(l)}$
in the $l$ variables $x_1, \dots, x_l$. Then the corresponding $\beta_l$
arising from inclusion/exclusion is ${{n}\choose{l}}\cdot A(l)$. Now we
want to show that under modest asymptotic conditions, the horrendous sum
in the definition of $A(l)$ can be replaced by the single term arising
from $w = (r,r,\dots, r)$, where $r=R-1$, exactly as in the preceding
argument.

\comment
The substitution $x_i \mapsto x$ induces an equality for each $s$, 
$$\sum_{w \in D(s)} \(\(1-lp + p \sum_{i=1}^l
x_i\)^{K}, \sum x^w \) = \((1 - lp + lp x \)
$$
\endcomment 

For $s < lr$ and $w = (a(j)) \in D(s)$, we bound
$$\eqalign{
\frac{ \(\(1-lp + p \sum_{i=1}^l x_i\)^{K}, x^w \)}{ \(\(1-lp + p
\sum_{i=1}^l x_i\)^{K}, \prod x_i^r \)} &=\frac{ {{K}\choose {a(1), a(2),
\dots, a(l)}} p^s (1-lp)^{K-s} }
{ {{K}\choose {r,r,\dots, r}} p^{lr} (1-lp)^{K-lr} }\cr
& = \frac{\frac{K!}{(K-s)!\prod_j a(j)!}}  {\frac{K!}{(K-lr)!(r!)^l} } \cdot \(\frac{1-lp}{p}\)^{lr-s} \cr
& = \frac{(r!)^l}{\prod_j a(j)!}\cdot \frac 1{\prod_{i=s+1}^{lr} (K-i)}
\cdot \(\frac{1-lp}{p}\)^{lr-s} \cr
& =  \frac{(r!)^l}{\prod_j a(j)!}\cdot  \(\frac{1-lp}{Kp}\)^{lr-s}
\exp\((l^2 r^2- s^2)/2K\)\cdot \( 1 + \Oh{r^3l^3/K^2}\).
}$$
With $ p = S/N$, $d = lr-s$ and $KS = NR\alpha(N)$, this translates to
(modulo a factor of $ 1 + \Oh{r^3l^3/K^2}$ (so we  assume $r(N) =
\oh{K(N)^{2/3}}$; generally, $l$ will be fixed),
$$\eqalign{
 \frac{(r!)^l}{\prod_j a(j)!}\cdot  \(\frac{1-lS/N}{KS/N}\)^{d} \exp\((l^2
r^2- s^2)/2K\)& \leq
 \(\frac{1-lS/N}{R\alpha(N)}\)^{d} \exp\((l^2 r^2- s^2)/2K\)
\frac{(r!)^l}{\prod_j a(j)!}\cr
& \leq  \(\frac{1}{R\alpha(N)}\)^{d} \frac{(r!)^l}{\prod_j a(j)!}.\cr
}$$
If we go back to the situation in (1) with $ m = l$ and $i_j = a(j)$, the term $\prod {S\choose {i_j}}$ expands as $S^s /\prod_j (a(j)!) \exp (-\sum a(j)^2/S)$ (up to multiplication by $1 \pm \oh 1$). Since $a(j) \leq R-1$ and $R^2 = \oh{S}$ ($S$ is now the number of bins of course), we see that $\sum_{s=0}^{rl-1} 1/$

We are thus in the same position as before, and we sum, first over $w \in
D(s)$, and then over $s <  lr$ (We always assume that  $N$ is chosen
sufficiently large that $l^2/K$ is miniscule.) This means we can replace
$\beta_l = {{n}\choose{l}}\cdot A(l)$ by
$$\eqalign{
\beta'_l & = {{N/S}\choose l}{{K} \choose{r,r,\dots, r}} p^{lr} (1-lp)^{K
-lr}\cr
& = \frac{(N/S)^{l}}{l!} \(1 + \Oh{l^2S/N} \) \frac{K!}{(r!)^l (K-lr)!}
p^{lr} (1-lp)^{K -lr}\cr
& = \frac1 {l!}\(\frac{N(Kp)^r}{r! S}\)^l \exp\(- lKp\) \(1 +\Oh{
\max\brcs{S/N, Kp^2}}\)\cr
&=  \frac1 {l!}\(\frac{N(Kp)^r}{r! S}\)^l \exp\(- lKp\) \(1 +\Oh{
\max\brcs{S/N, Kp^2}}\)\cr
& = \frac{\(\frac{N(\alpha(N)R)^{R-1})}{S(R-1)!}\)^l}{l!}  \(1 +\Oh{
\max\brcs{S/N, Kp^2}}\).\cr
}$$
This is exactly the same as xxx (assuming $\alpha(N)RS = \oh{N}$). Hence
the asymptotic solutions are identical, including  the asymmetric phase
transitions.

Statement of results for balls and bins.

Now to do the balls and bins maximal problem. It's different!!

We can also consider the  maximum problem, and run into results that are
already known in the balls and bins situation. Given a subdivision of
$\brcs{0,1,2, \dots, N}$ into the $\flo{N/S}$ consecutive intervals $I_k=
[kS, (k+1)S -1]$ ($k = 0,1, \dots, \flo {N/S}-1$, ignoring the last short
interval if $S$ does not divide $N$), what is the likelihood that a set of
cardinality $K$ hits at least one $I_k$ in at least $A+1$ points? Here $K$
is viewed as a parameter, and obviously we may as well insist that $A+1
\leq K\leq NA/S$ (so that, unlike the case involving the minimum
occupancy, $K$ is going to be $\oh{N}$ (assuming $S(N) \to \infty$, which
we do). As before, $A, S, K$ are functions of $N$, and we will
(eventually) let $N$ wander off to infinity.

This is easy to quantify, but not so easy to solve. Abbreviate $l =
\flo{N/S}$, the number of intervals (in analogy with the balls and bins
problem, this should be $m$, but we've already used that in the
$\beta_m$s). The negation is that the set $U$ satisfies $|I_k \cap U| \leq
A$ for all $k = 0,1,\dots, \flo{N/S}-1$. Let $D_A(s) = \sum_{(a(j)) \in
\Z^l_+}{a(j) \leq A \text{ for all $j$ and }\sum_{j=1}^l a(j)  = s}$.  The
probability of the negation is
$$\eqalign{
\sum_{\cup_{s=0}^K D(s) }\frac{\prod_{j=1}^l {{S} \choose a(j)}}{N\choose
K}\cr
\& = \frac{\sum_{s= 0}^{K} \sum_{(a(j)) \in D(s)} \frac{\prod_{j=1}^l {{S}
\choose a(j)}}{N\choose K}\cr
& = \frac{\sum_{s = 0}^K \(\( \sum_{i=0}^A x^i {S \choose i} \)^l, x^s)}{N
\choose K}\cr
}$$
Since the polynomial $f_A = \sum_0^A x^i {S \choose i}$ is strongly
unimodal (that is, the distribution of its coefficients is log concave;
that is so because $f_A$ is just an initial segment of $(1+x)^{S}$), so is
its $l$th power, and therefore the sequence $(h(0), h(1), \dots, h(K))$
given by
$$
h(s) =  \sum_{(a(j)) \in D(s)} \frac{\prod_{j=1}^l {{S} \choose a(j)}
$$
is log concave. If $(b(i))_{i=0}^M$ is a distribution of nonnegative real
number, and $f = \sum b_i x^i$, then the mean of the distribution is
$(\sum b_i i) = f'(1)/f(1)$, denoted $\mu (f)$. So the mean of the
distribution for $h$f_A'(1)/f_A(1)$) is $l\cdot \mu(f) = l \cdot \sum_0^A
{S\choose i}i/ \sum_{S \choose i}$. Using the trivial result (xxx), by
computing the ratio of the largest term to the second largest (assuming $A
< S/2)$,  we have
$$\eqalign{
{S \choose A }(1 + A/(S-A))& \leq \sum_0^A {S \choose i} \leq {S \choose A
}(1 + (A/(S-2A))\cr
{S \choose A }A(1 + (A-1)/(S-A))& \leq \sum_0^A{S \choose i} i\leq {S
\choose A }A(1 + (A-1)(S-2A+1).\cr
}$$
We assume that $A = \oh{S}$ (we may need $A = \oh{\sqrt S}$????), so that
$\mu(f_A)$ is very close to $A$ (it must be strictly less than $A$, but
asymptotically, at least if $A = \oh{S}$, it tends to $A$). Hence
$$
\(1 - \frac{A(S+A-2)}{(S-A)^2}\)A=\(\frac{1+ \frac{A-1}{S-A}}{1 +
\frac{A}{S-2A}}\)A\geq \mu(f_A) \leq A \(\frac{1+ \frac {A-1}{S-A}{1+
\frac A{S-A}}}\)= A \( 1 - \frac{A^2 + 1}{S(S-2A+1)}\).
$$
(we can easily make these inequalities sharper by considering the third
highest terms, etc; but there is no need to go too much farther). Thus the
mean of the distribution associated to $h$, that is, $\mu(f_A)^l)$ is at
least $(N/S)(A-1)$. Unless $K \geq (N/S)(A-1)$ (in which case it is almost
certain that the corresponding set hit at least one of the $I_k$ in more
than $A$ points), $K$ is much smaller than the mean. By Darroch's theorem,
the mean and mode are within one of each other, so that $h$ is increasing
on the interval $\brcs{0,1,\dots, K}$.

The variance, $\sigma^2:= V(f) = f'(1)/f(1) - f'(1)^2/f(1)^2 +
f''(1)/f(1)$, is also easily computed for the distribution of $h$. Again
using log concavity and computing the ratio of the largest (= last) term
to the second largest (= second last), we see that ${S \choose
A}^{-2}((f_A''(1) + f_A'(1))f(1) - f_A'(1)^2)\geq (AS^2 - 2A^4S + 3A^3 S
\pm\oh{A^3S})/(S-2A)(S-2A+ 1)(S-A+1)$, essentially nothing.
$$
{S \choose A }A^2(1+(A-1)^2/A(S-A+1)) \sum_0^A{S \choose i} i^2 \leq {S
\choose A }A^2(1 + (A-1)^2/(AS - 2A^2 + 3A + 1)).
$$
Seemingly we obtain $A^4/S^2-A/(S-A) <\sigma^2 <  A^4/(S-A)^2$ or
something like that. If $A = \oh{S^2}$, $\sigma$ is essentially zero. On
the other hand, if not, then we can use Gauss inequality (for unimodal
distributions), expressed in the form $\text{Pr}(|X -\mu| > k\sigma) \leq
4/9k^2$. Here $\mu $ is roughly $NA/S(1-A/(A-S))$, $\sigma \sim NA^2/S^2$
(we assume at least $A = \oh{S}$), and on setting $NA/S-K \sim \sigma k$,
we deduce $k \sim S^2(NA/S - K)/NA^2 = (NAS- S^2K)/NA^2 $. Set $t = K/A$
(so $1 \leq t$); then $k \sim (S/A)(1- St/N)$. Then the probability is at
most $(A^2/S^2) (1-St/N)$. So the probability of failure (remember, this
was the negation of the original question) is asymptotically zero

Something is wrong here. As $K$ increases, probability should go to zero.
Hence if $(K-A)/A $ is not $\Oh{A/S}$, the generic  set of size $K$ will
hit at least one of the $I_k$ i